\documentclass[onecolumn,11pt]{article}
\usepackage[top=1in, bottom=1in, left=1in, right=1in]{geometry}
\usepackage{amsmath,amsfonts,amscd,amssymb}
\usepackage{graphicx}
\usepackage{epstopdf}
\usepackage{overpic}
\usepackage{cancel}
\usepackage{rotating}
\usepackage{url}
\usepackage{caption}
\usepackage{color}
\usepackage{rotating}
\usepackage{multirow}
\usepackage{wrapfig}
\usepackage{mathtools}
\usepackage{subeqnarray}
\usepackage{setspace}
\usepackage{palatino} 
\usepackage[bottom,flushmargin,hang,multiple]{footmisc}
\usepackage{lipsum}
\newcommand\blfootnote[1]{%
  \begingroup
  \renewcommand\thefootnote{}\footnote{#1}%
  \addtocounter{footnote}{-1}%
  \endgroup
}

\definecolor{header1}{cmyk}{0,0,0,1}
\newcommand{\HLtext}{\sffamily\bfseries}

\newcommand{\bm}{\mathbf{m}}

\title{\LARGE{\vspace{-.5in}\textbf{Discovery of Nonlinear Multiscale Systems:  Sampling Strategies and Embeddings}}\vspace{-.15in}}
\author{\normalsize{Kathleen P. Champion$^{1*}$, Steven L. Brunton$^2$, J. Nathan Kutz$^1$}\\
\footnotesize{$^1$ Department of Applied Mathematics, University of Washington, Seattle, WA 98195, United States}\\
\footnotesize{$^2$ Department of Mechanical Engineering, University of Washington, Seattle, WA 98195, United States\vspace{-.1in}}
}
\date{}
\begin{document}
\maketitle

\blfootnote{$^*$ Corresponding author (kpchamp@uw.edu).\\ \noindent \textbf{Code:}  https://github.com/kpchamp/MultiscaleDiscovery}
\vspace{-.2in}
\begin{abstract}
A major challenge in the study of dynamical systems is that of model discovery: turning data into models that are not just predictive, but provide insight into the nature of the underlying dynamical system that generated the data. This problem is made more difficult by the fact that many systems of interest exhibit diverse behaviors across multiple time scales. We introduce a number of data-driven strategies for discovering nonlinear multiscale dynamical systems and their embeddings from data. We consider two canonical cases: (i) systems for which we have full measurements of the governing variables, and (ii) systems for which we have incomplete measurements. For systems with full state measurements, we show that the recent sparse identification of nonlinear dynamical systems (SINDy) method can discover governing equations with relatively little data and introduce a sampling method that allows SINDy to scale efficiently to problems with multiple time scales.  Specifically, we can discover distinct governing equations at slow and fast scales. For systems with incomplete observations, we show that the Hankel alternative view of Koopman (HAVOK) method, based on time-delay embedding coordinates, can be used to obtain a linear model and Koopman invariant measurement system that nearly perfectly captures the dynamics of nonlinear quasiperiodic systems. We introduce two strategies for using HAVOK on systems with multiple time scales. Together, our approaches provide a suite of mathematical strategies for reducing the data required to discover and model nonlinear multiscale systems.\\

\noindent\emph{Keywords--}
model discovery, dynamical systems, Hankel matrix, sparse regression, time-delay embedding, Koopman theory, multiscale dynamics
\end{abstract}


 \section{Introduction}
With recent advances in sensor technologies and data collection techniques in fields such as biology, neuroscience, and engineering, there are increasingly many systems for which we have significant quantities of measurement data but do not know the underlying governing equations. Thus a major challenge is that of model discovery: finding physical models, governing equations and/or conserved quantities that drive the dynamics of a system. Learning governing models from data has the potential to dramatically improve our understanding of such systems, as well as our ability to predict and control their behavior. However many systems of interest are highly complex, exhibiting nonlinear dynamics and multiscale phenomena, making the problem of model discovery especially difficult.   To this end, we develop a diverse set of mathematical strategies for discovering nonlinear multiscale dynamical systems and their embeddings directly from data, showing their efficacy on a number of example problems.

The recent and rapid increase in the availability of measurement data has spurred the development of many data-driven methods for modeling and predicting dynamics. At the forefront of data-driven methods are deep neural networks (DNNs).  DNNs not only achieve superior performance for tasks such as image classification \cite{krizhevsky_imagenet_2012}, but they have also been shown to be effective for future state prediction of dynamical systems~\cite{li_extended_2017,vlachas_data-driven_2018,yeung_learning_2017,Takeishi2017nips,Wehmeyer2017arxiv,Mardt2017arxiv,lusch_deep_2017,raissi2018multistep}.   A key limitation of DNNs, and similar data-driven methods, is the lack of interpretability of the resulting model:  they are focused on prediction and do not provide governing equations or clearly interpretable models in terms of the original variable set.  An alternative data-driven approach uses symbolic regression to identify directly the structure of a nonlinear dynamical system from data \cite{bongard_automated_2007,schmidt_distilling_2009,cornforth_symbolic_2012}. This works remarkably well for discovering interpretable physical models, but symbolic regression is computationally expensive and can be difficult to scale to large problems.  However, the discovery process can be reformulated in terms of sparse regression~\cite{brunton_discovering_2016}, providing a computationally tractable alternative.

Regardless of the data-driven strategy used, many complex systems of interest exhibit behavior across multiple time scales, which poses unique challenges for modeling and predicting their behavior. It is often the case that while we are primarily interested in macroscale phenomena, the microscale dynamics must also be modeled and understood, as they play a role in driving the macroscale behavior.  
The macroscale dynamics in turn drive the microscale dynamics, thus producing a coupling whereby the dynamics at different time scales feedback into each other. This can make dealing with multiscale systems particularly difficult unless the time scales are disambiguated in a principled way. There is a significant body of research focused on modeling multiscale systems: notably the heterogeneous multiscale modeling (HMM)  framework and equation-free methods for linking scales \cite{kevrekidis_equation-free_2003,weinan_heterognous_2003,weinan_principles_2011}. Additional work has focused on testing for the presence of multiscale dynamics so that analyzing and simulating multiscale systems is more computationally efficient \cite{froyland_computational_2014,froyland_trajectory-free_2016}. Many of the same issues that make modeling multiscale systems difficult can also present challenges for model discovery and system identification. This motivates the development of specialized methods for performing model discovery on problems with multiple time scales, taking into account the unique properties of multiscale systems.

In this paper, we consider strategies for the discovery of nonlinear multiscale dynamical systems and their embeddings. In particular, we consider two cases: systems for which the full state measurements are available, and those for which we have incomplete measurements and latent variables. In Section~\ref{sec:full_state_measurements}, we consider the case where we have full state measurements for the system of interest. We focus on the recent sparse identification of nonlinear dynamical systems (SINDy) algorithm \cite{brunton_discovering_2016} and establish a baseline understanding of how much data the method requires. We then introduce a sampling method that enables the algorithm to scale efficiently for multiscale problems, thus disambiguating between the distinct dynamics that occur at fast and slow scales. Our method allows SINDy to perform well with approximately the same amount of training data, even as the fast and slow time scales diverge. In Section~\ref{sec:incomplete_measurements}, we consider systems for which we do not have full state measurements. We show that for quasiperiodic systems, the Hankel alternative view of Koopman (HAVOK) method \cite{brunton_chaos_2017} provides a closed linear model using time-delay embeddings that allows for highly accurate reconstruction and long-term prediction of the dynamics. These models do not need the nonlinear forcing term that was necessary for the chaotic systems studied in \cite{brunton_chaos_2017}. We then provide two strategies for using HAVOK to model multiscale systems. The first is a strategy for constructing our data matrices so that we can create an accurate model for the system with less data. The second is an iterative method for separating fast and slow time scales. Together, the approaches in this paper enable a significant reduction in the data requirement for the discovery and modeling of nonlinear multiscale dynamical systems.

\begin{figure}[t]
  \vspace{.1in}
  \centering
  \begin{overpic}[width=.98\linewidth]{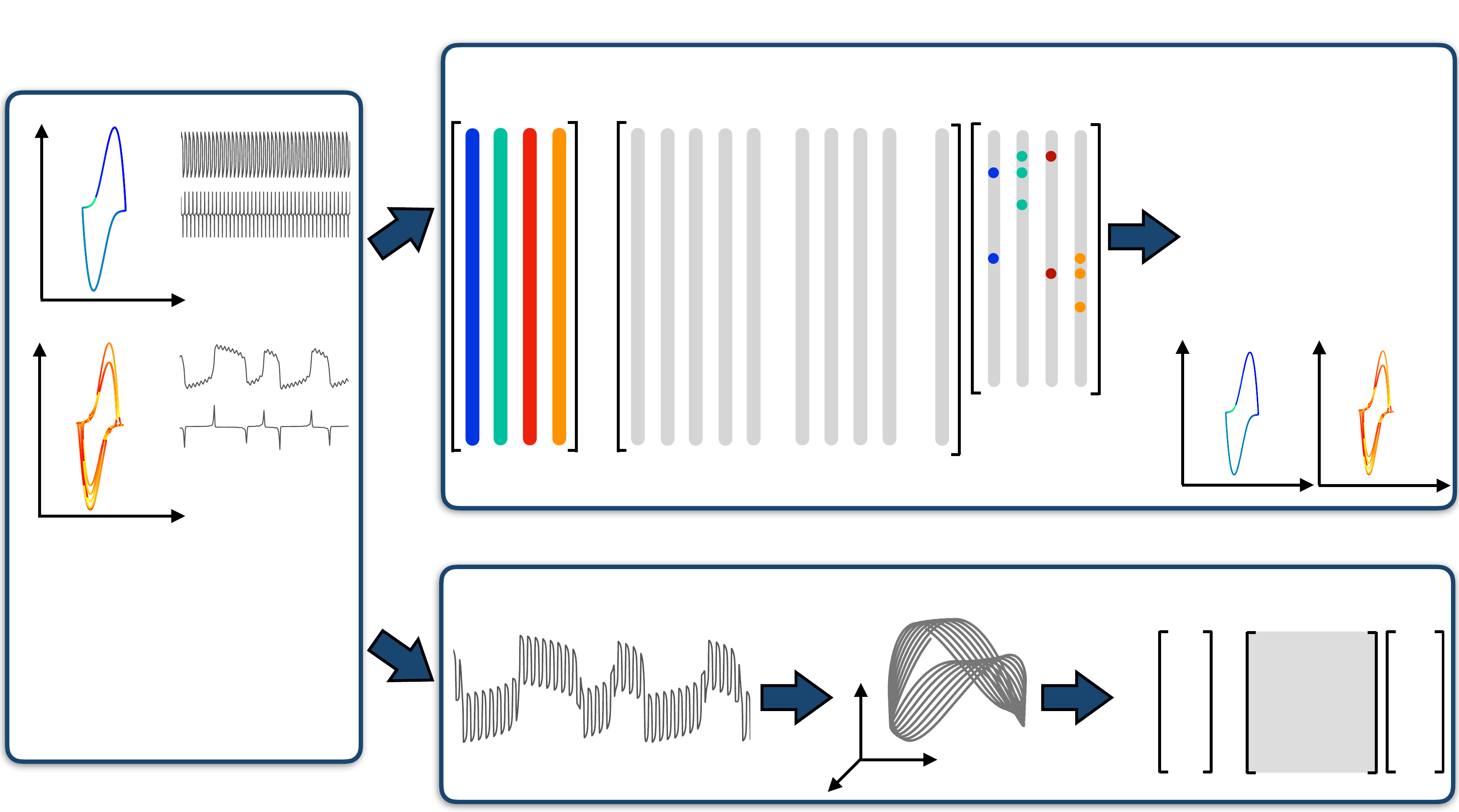}
  \put(3,51){\color{header1}\HLtext (a) True System}
  \put(1,41){\tiny $y_1$}
  \put(1,26){\tiny $y_2$}
  \put(7,34){\tiny $x_1$}
  \put(7,19.5){\tiny $x_2$}
  \put(0,11){\small \parbox{0.25\linewidth}{,
      \begin{align*}
    \dot{x}_1 &= y_1 + c_1 x_2 \\
    \dot{y}_1 &= \mu(1-x_1^2)y_1 - x_1 \\
    \dot{x}_2 &= y_2 + c_2 x_1 \\
    \dot{y}_2 &= \mu(1-x_2^2)y_2 - x_2
    \end{align*}}}
  \put(50,54){\color{header1}\HLtext (b) Full state measurements}
  \put(32,50.2){\small SINDy}
  \put(81,50.2){\small Identified system}
  \put(82,41){\small \parbox{0.15\linewidth}{
      \begin{align*}
    \dot{x}_1 &= \mathbf{\Theta}(\mathbf{x}^T)\xi_1 \\
    \dot{y}_1 &= \mathbf{\Theta}(\mathbf{x}^T)\xi_2 \\
    \dot{x}_2 &= \mathbf{\Theta}(\mathbf{x}^T)\xi_3 \\
    \dot{y}_2 &=\mathbf{\Theta}(\mathbf{x}^T)\xi_4
    \end{align*}}}
  \put(32,48){\tiny $\dot{x}_1$}
  \put(34,48){\tiny $\dot{y}_1$}
  \put(35.7,48){\tiny $\dot{x}_2$}
  \put(37.5,48){\tiny $\dot{y}_2$}
  \put(43.2,48){\tiny $1$}
  \put(44.5,48){\tiny $x_1$}
  \put(46.5,48){\tiny $y_1$}
  \put(48.5,48){\tiny $x_1^2$}
  \put(50.3,48){\tiny $x_1y_1$}
  \put(54.2,48){\tiny $x_2$}
  \put(56,48){\tiny $y_2$}
  \put(57.8,48){\tiny $x_2^2$}
  \put(59.7,48){\tiny $x_2y_2$}
  \put(63.2,48){\tiny $x_1y_2^2$}
  \put(67.7,48){\tiny $\xi_1$}
  \put(69.6,48){\tiny $\xi_2$}
  \put(71.5,48){\tiny $\xi_3$}
  \put(73.4,48){\tiny $\xi_4$}
  \put(40.1,35){$=$}
  \put(52.2,35){\footnotesize $\cdots$}
  \put(61.7,35){\footnotesize $\cdots$}
  \put(34,22.5){$\mathbf{X}$}
  \put(51,22.5){$\mathbf{\Theta}(\mathbf{X})$}
  \put(70.5,26){$\mathbf{\Xi}$}
  \put(36,18){\color{header1}\HLtext (c) Incomplete measurements and latent variables}
  \put(31,14.5){\small Measure $x_1 + x_2$}
  \put(58,14.5){\small Delay embed}
  \put(78,14.5){\small Regression model}
  \put(76.9,7.5){$\frac{\mathrm{d}}{\mathrm{d}t}$}
  \put(79.5,7.5){\scriptsize $\begin{array}{c} v_1 \\ v_2 \\ \vdots \\ v_r \end{array}$}
  \put(83.5,7.5){\footnotesize $=$}
  \put(88.5,7){\Large $\mathbf{A}$}
  \put(95.5,7.5){\scriptsize $\begin{array}{c} v_1 \\ v_2 \\ \vdots \\ v_r \end{array}$}
  \end{overpic}
  \caption{Schematic showing methods for modeling dynamical systems with multiple time scales. (a) An example system of two coupled Van der Pol oscillators, one fast and one slow. We show the phase portraits, time dynamics, and full set of underlying equations. (b) When full state measurements are available, sparse identification of nonlinear dynamics (SINDy) can be used to discover the system. SINDy uses sparse regression to identify governing equations. (c) When we have incomplete measurements or latent variables driving the dynamics, alternative methods such as Hankel alternative view of Koopman (HAVOK) must be used. HAVOK uses time-delay coordinates to produce a delay embedded attractor. Regression is then used to produce a linear model for the dynamics on the attractor.}
  \label{fig:1}
\end{figure}
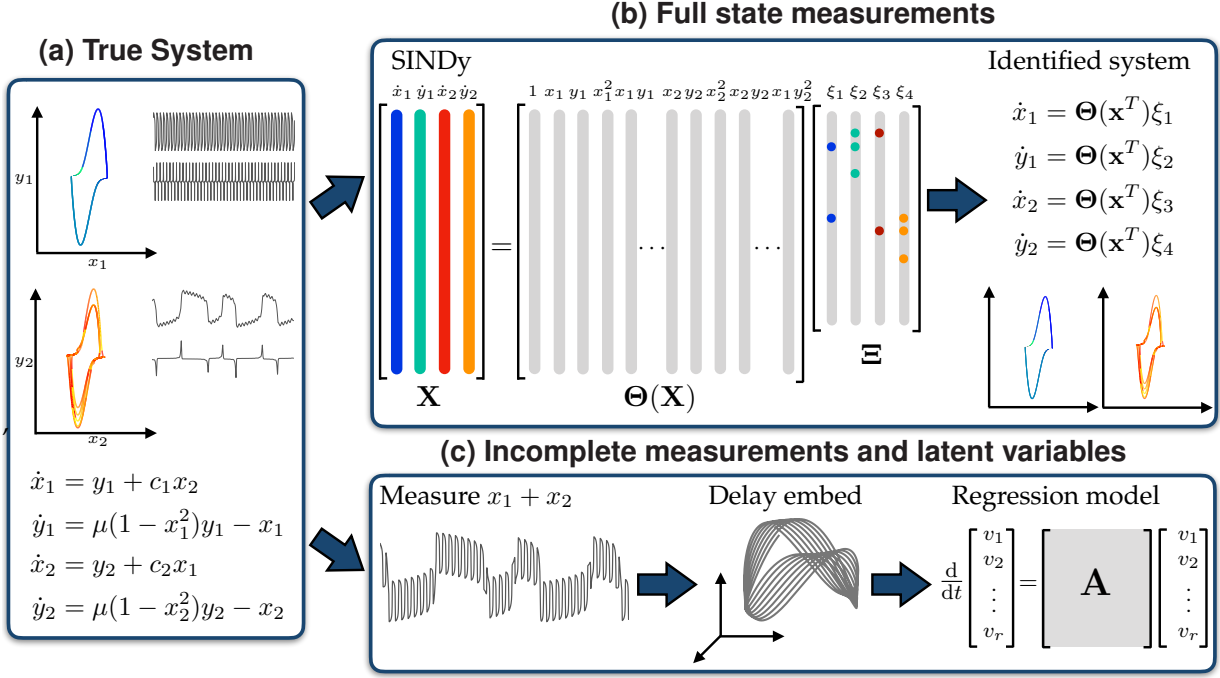

\section{Systems with full state measurements}
\label{sec:full_state_measurements}
We consider dynamical systems of the form
\begin{equation}
  \frac{d}{dt}\mathbf{x}(t) = \mathbf{f}(\mathbf{x}(t))
  \label{eq:dynamical_system}.
\end{equation}
Here $\mathbf{x}(t) \in \mathbb{R}^n$ is the state of the system at time $t$ and $\mathbf{f}$ defines the dynamics of the system.  A key challenge in the study of dynamical systems is finding governing equations; that is, to find the form of $\mathbf{f}$ given measurements of $\mathbf{x}(t)$.   

In this section we investigate the sparse identification of nonlinear dynamical systems (SINDy) algorithm, a recent method for discovering the governing equations of nonlinear dynamical systems of the form (\ref{eq:dynamical_system}) from measurement data \cite{brunton_discovering_2016}. The SINDy algorithm relies on the assumption that $\mathbf{f}$ has only a few active terms; it is sparse in the space of all possible functions of $\mathbf{x}(t)$. Given snapshot data
\begin{equation}
  \mathbf{X} = \left(\begin{array}{cccc}
    x_1(t_1) & x_2(t_1) & \cdots & x_n(t_1) \\
    x_1(t_2) & x_2(t_2) & \cdots & x_n(t_2) \\
    \vdots & \vdots & \ddots & \vdots \\
    x_1(t_m) & x_2(t_m) & \cdots & x_n(t_m) \\
  \end{array}\right), \quad
  \dot{\mathbf{X}} = \left(\begin{array}{cccc}
    \dot{x}_1(t_1) & \dot{x}_2(t_1) & \cdots & \dot{x}_n(t_1) \\
    \dot{x}_1(t_2) & \dot{x}_2(t_2) & \cdots & \dot{x}_n(t_2) \\
    \vdots & \vdots & \ddots & \vdots \\
    \dot{x}_1(t_m) & \dot{x}_2(t_m) & \cdots & \dot{x}_n(t_m) \\
  \end{array}\right), \nonumber
\end{equation}
the SINDy method builds a library of candidate nonlinear functions $\mathbf{\Theta}(\mathbf{X}) = [\theta_1(\mathbf{X}) \cdots \theta_p(\mathbf{X})]$. There is great freedom in choosing the candidate functions, although they should include the terms that make up the function $\mathbf{f}(\mathbf{x})$ in (\ref{eq:dynamical_system}).    Polynomials are of particular interest as they are often key elements in canonical models of dynamical systems.  Another interpretation is that the SINDy algorithm discovers the dominant balance dynamics of the measured system, which is often of a polynomial form due to a Taylor expansion of a complicated nonlinear function.
The algorithm then uses thresholded least squares to find sparse coefficient vectors $\mathbf{\Xi} = (\mathbf{\xi}_1\ \mathbf{\xi}_2\ \cdots\ \mathbf{\xi}_n )$ to approximately solve
\begin{equation}
  \dot{\mathbf{X}} = \mathbf{\Theta}(\mathbf{X})\mathbf{\Xi}. \nonumber
\end{equation}
The SINDy algorithm is capable of identifying the nonlinear governing equations of the dynamical system, provided that we have measured the correct states $\mathbf{x}(t)$ that contribute to the dynamics and the library $\mathbf{\Theta}$ is rich enough to span the function $\mathbf{f}$. We address systems for which we do not have full state measurements in Section~\ref{sec:incomplete_measurements}. In this section we focus on systems for which we have full state measurements and explore the performance of the SINDy algorithm on both uniscale and multiscale systems.

\subsection{Uniscale dynamical systems}
\label{sec:uniscale1}

Before considering multiscale systems, we establish a baseline understanding of the data requirements of the SINDy algorithm. While previous results have assessed the performance of SINDy in various settings \cite{mangan_model_2017,kaiser_sparse_2017}, so far none have looked explicitly at how much data is necessary to correctly identify a system. We determine the data requirements of SINDy on four example systems: the periodic Duffing and Van der Pol oscillators, and the chaotic Lorenz and Rossler systems.  
To assess how quickly SINDy can correctly identify a given system, we look at its performance on measurement data from a single trajectory. Each example system has dynamics that evolve on an attractor. We choose data from time points after the system has converged to the attractor and quantify the sampling rate and duration in relation to the typical time $T$ it takes the system to make a trip around some portion of the attractor. We refer to this as a ``period'' of oscillation. This gives us a standard with which to compare the data requirements among multiple systems.
We show that in all four models considered, SINDy can very rapidly discover the underlying dynamics of the system, even if sampling only a fraction of a period of the attractor or oscillation.

Our procedure is as follows. For all systems, we simulate a single trajectory at a sampling rate $r=2^{18}$ samples/period, giving us a time step $\Delta t=T/r$ (note $T$ is defined differently for each system). We then subsample the data at several rates $r_\text{sample} = 2^5,2^6,\dots,2^{18}$ and durations on different portions of the attractor. For each set of subsampled data, we train a SINDy model and determine if the identified model has the correct set of nonzero coefficients, meaning SINDy has identified the proper form of the dynamics. We are specifically looking for the correct structural identification of the coefficient matrix $\mathbf{\Xi}$, rather than precise coefficient values. In general we find that if SINDy correctly identifies the active coefficients in the model, the coefficients are reasonably close to the true values. Our library of SINDy candidate functions includes polynomial terms of up to order $3$. In this section, we choose the coefficient threshold for the iterative least squares algorithm to be $\lambda=0.1$. The results indicate the sampling rate and length of time one must sample for SINDy to correctly discover a dynamical system.

When assessing the performance of SINDy, it is important to consider how measurements of the derivative are obtained. While in some cases we may have measurements of the derivatives, most often these must be estimated from measurements of $\mathbf{x}$. In this work, we consider the low noise case and are thus able to use a standard center difference method to estimate derivatives. With noisy data, more sophisticated methods such as the total variation regularized derivative may be necessary to obtain more accurate estimates of the derivative \cite{chartrand_numerical_2011,brunton_discovering_2016}.  Indeed, accurate computations of the derivative are critical to the success of SINDy.

To assess performance, we consider the duration of sampling required to identify the system at each chosen sampling rate. The exact duration required depends on which portion of the attractor the trajectory is sampled from; thus we take data from different portions of the attractor and compute an average required duration for each sampling rate. 
This provides insight into how long and at what rate we need to sample in order to correctly identify the form of the dynamics. Surprisingly, in all four example systems, {\em data from less than a full trip around the attractor is typically sufficient for SINDy to identify the correct form of the dynamical system}.

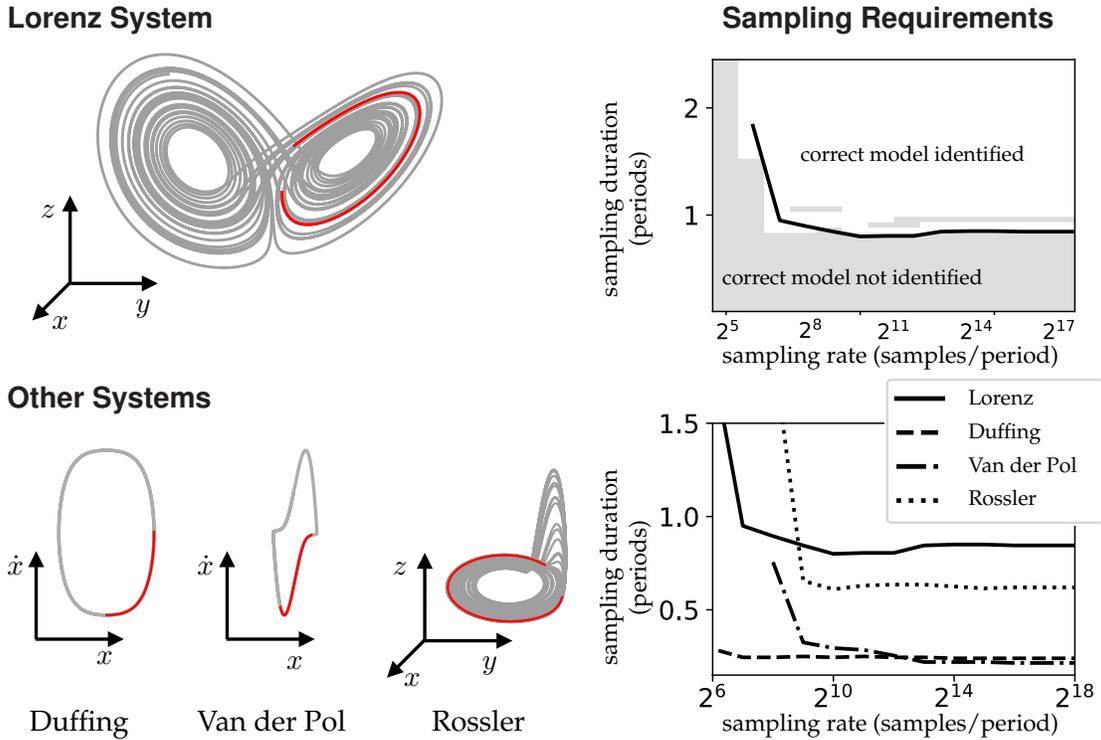
\begin{figure}[t]
  \vspace{.1in}
  \centering
  \begin{overpic}[width=.9\linewidth]{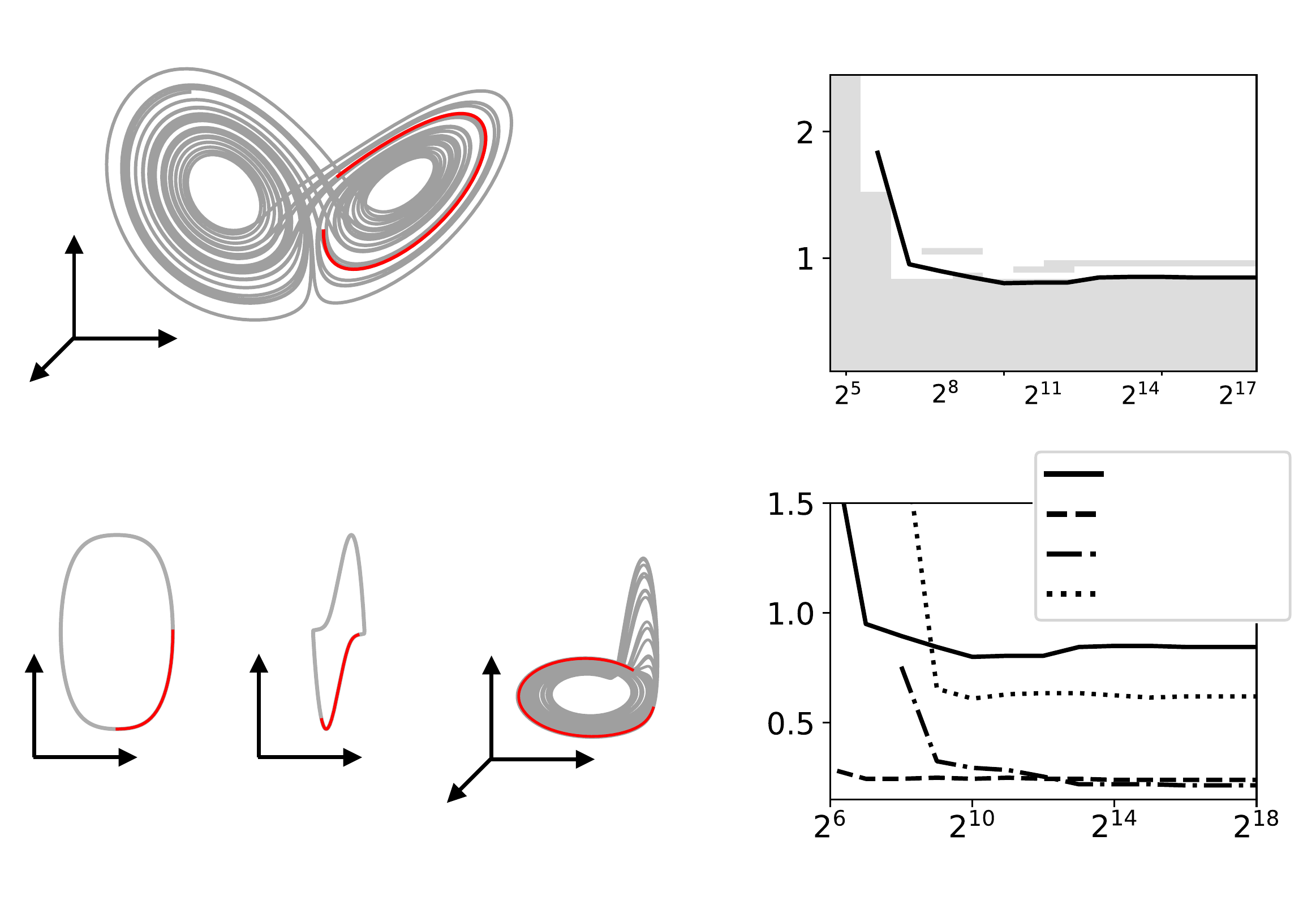}
  \put(0,66){\color{header1}\HLtext Lorenz System}
  \put(64,66){\color{header1}\HLtext Sampling Requirements}
  \put(4,39){$x$}
  \put(11.5,40.5){$y$}
  \put(3,49){$z$}
  \put(64,43){\scriptsize correct model not identified}
  \put(71,54){\scriptsize correct model identified}
  \put(53,40){\footnotesize \rotatebox{90}{\parbox{10em}{\centering sampling duration (periods)}}}
  \put(64,36){\footnotesize sampling rate (samples/period)}
  \put(0,32){\color{header1}\HLtext Other Systems}
  \put(2,3){Duffing}
  \put(8,9){$x$}
  \put(0,17){$\dot{x}$}
  \put(17,3){Van der Pol}
  \put(25,8.5){$x$}
  \put(17,17){$\dot{x}$}
  \put(38,3){Rossler}
  \put(35.5,7){$x$}
  \put(42.5,8.5){$y$}
  \put(34.5,17){$z$}
  \put(86,32.2){\scriptsize Lorenz}
  \put(86,29.2){\scriptsize Duffing}
  \put(86,26.2){\scriptsize Van der Pol}
  \put(86,23.2){\scriptsize Rossler}
  \put(64,3){\footnotesize sampling rate (samples/period)}
  \put(53,7){\footnotesize \rotatebox{90}{\parbox{10em}{\centering sampling duration (periods)}}}
  \end{overpic}
  \caption{Data requirements of SINDy for four example systems: the Lorenz system, Duffing oscillator, Van der Pol oscillator, and Rossler system. We plot the attractor of each system in gray. The portion of the attractor highlighted in red indicates the average portion we must sample from to discover the system using SINDy. For all systems, we see that we do not need to sample from the full attractor. Plots on the right indicate how the sampling rate affects the duration we must sample to obtain the correct model. In the bottom plot, we show for each system the average sampling duration necessary to identify the correct model at various sampling rates. The top plot provides a detailed look at SINDy's performance on the Lorenz system. (top right) Gray shading indicates regions where the correct model was not identified. The black curve indicates on average the sampling duration necessary at each sampling rate (and is the same curve plotted for the Lorenz system in the bottom plot).}
  \label{fig:2}
\end{figure}

\subsubsection{Lorenz system}
\label{sec:lorenz}
As a first example, consider the chaotic Lorenz system:
\begin{subequations}
\begin{align}
\dot{x} &= \sigma(y-x) \label{eq:lorenz} \\
\dot{y} &= x(\rho - z) - y  \\
\dot{z} &= xy - \beta z, 
\end{align}
\end{subequations}
with parameters $\sigma=10$, $\rho=28$, and $\beta=8/3$. The dynamics of this system evolve on an attractor shown in Figure~\ref{fig:2}. At the chosen parameter values the Lorenz attractor has two lobes, and the system can make anywhere from one to many cycles around the current lobe before switching to the other lobe. For the Lorenz system, we define the period $T$ to be the typical time it takes the system to travel around one lobe. At the chosen parameter values, we determine that $T \approx 0.759$ (calculated by averaging over many trips around the attractor); however, because of the chaotic nature of the system, the time it takes to travel around one lobe is highly variable.

In the top right panel of Figure~\ref{fig:2}, we plot the sampling rates and durations at which the system is correctly identified. The black curve shows the average duration of recorded data necessary to discover the correct model at each sampling rate. As the sampling rate increases, SINDy is able to correctly identify the dynamics with a shorter recording duration. In our analysis, we find that the system has both a baseline required sampling rate and a baseline required duration: if the sampling rate is below the baseline, increasing the duration further does not help discover the model. Similarly, if the duration is below the baseline, increasing the sampling rate does not help discover the model. For our chosen parameter values, we find the average baseline sampling duration to be 85\% of a period. Depending on the portion of the attractor sampled from, this duration ranges from 70-110\% of a period. The portion of the attractor covered by this average baseline duration is shown highlighted in red on the Lorenz attractor in Figure~\ref{fig:2}.

\subsubsection{Duffing oscillator}
Next we consider the Duffing oscillator, which can be written as a two-dimensional dynamical system:
\begin{subequations}
\begin{align}
\dot{x} &= y  \\
\dot{y} &= -\delta y - \alpha x - \beta x^3 .
\end{align}
\end{subequations}
We consider the undamped equation with $\delta=0$, $\alpha=1$, and $\beta=4$. The undamped Duffing oscillator exhibits regular periodic oscillations. The parameter $\beta > 0$ controls how nonlinear the oscillations are. At the selected parameter values, the system has period $T \approx 3.179$. We simulate the Duffing system, then apply the SINDy algorithm to batches of data subsampled at various rates and durations on different portions of the attractor.

The phase portrait and sampling requirements for the Duffing oscillator are shown in Figure~\ref{fig:2}. With a sufficiently high sampling rate, SINDy requires approximately 25\% of a period on average in order to correctly identify the dynamics. This portion of the attractor is highlighted in red on the phase portrait in Figure~\ref{fig:2}. Depending on where on the attractor we sample from, the required duration ranges from about 15-35\% of a period.

\subsubsection{Van der Pol oscillator}
\label{sec:vdp}
We next look at the Van der Pol oscillator, given by
\begin{subequations}
\begin{align}
\dot{x}_1 &= x_2 \label{eq:vdp} \\
\dot{x}_2 &= \mu (1-x_1^2)x_2 - x_1. 
\end{align}
\end{subequations}
The parameter $\mu>0$ controls the degree of nonlinearity; we use $\mu = 5$. At this parameter value we have period $T \approx 11.45$. We simulate the Van der Pol oscillator and again apply the SINDy method to sets of subsampled data.

The phase portrait and average sampling requirements for the Van der Pol oscillator are shown in Figure~\ref{fig:2}. The average baseline duration is around 20\% of a period. For this system in particular, the baseline duration is highly dependent on where the data is sampled from. If samples are taken during the fast part of the oscillation, the system can be identified with as little as 5\% of a period; sampling during the slower part of the oscillation requires as much as 35\% a period.

\subsubsection{Rossler system}
As a final example, consider the Rossler system
\begin{subequations}
\begin{align}
\tau\dot{x} &= -y-z \label{eq:rossler} \\
\tau\dot{y} &= x + ay  \\
\tau\dot{z} &= b + z(x-c), 
\end{align}
\end{subequations}
with parameters $a=0.1$, $b=0.1$, and $c=14$. Note we include a time constant $\tau$, with value $\tau=0.1$. The dynamics of this system evolve on the chaotic attractor shown in Figure~\ref{fig:2}. We define the period of oscillation in this system as the typical time it takes to make one trip around the portion of the attractor in the $x-y$ plane, which at these parameter values (with $\tau=0.1$) is $T \approx 6.14$. We simulate the system and apply the SINDy method as in previous examples.

The sampling requirements for the Rossler system are shown in Figure~\ref{fig:2}. Depending on the portion of the attractor we sample from, the baseline duration ranges from 35-95\% of a period of oscillation. Remarkably, SINDy can identify the system without any data from when the system leaves the $x-y$ plane. The average baseline sampling duration (65\% of a period) is highlighted in red on the Rossler attractor in Figure~\ref{fig:2}.

\subsection{Multiscale systems}
\label{sec:multiscale1}

\begin{figure}[t]
\vspace{.2in}
  \centering
  \begin{overpic}[width=\linewidth]{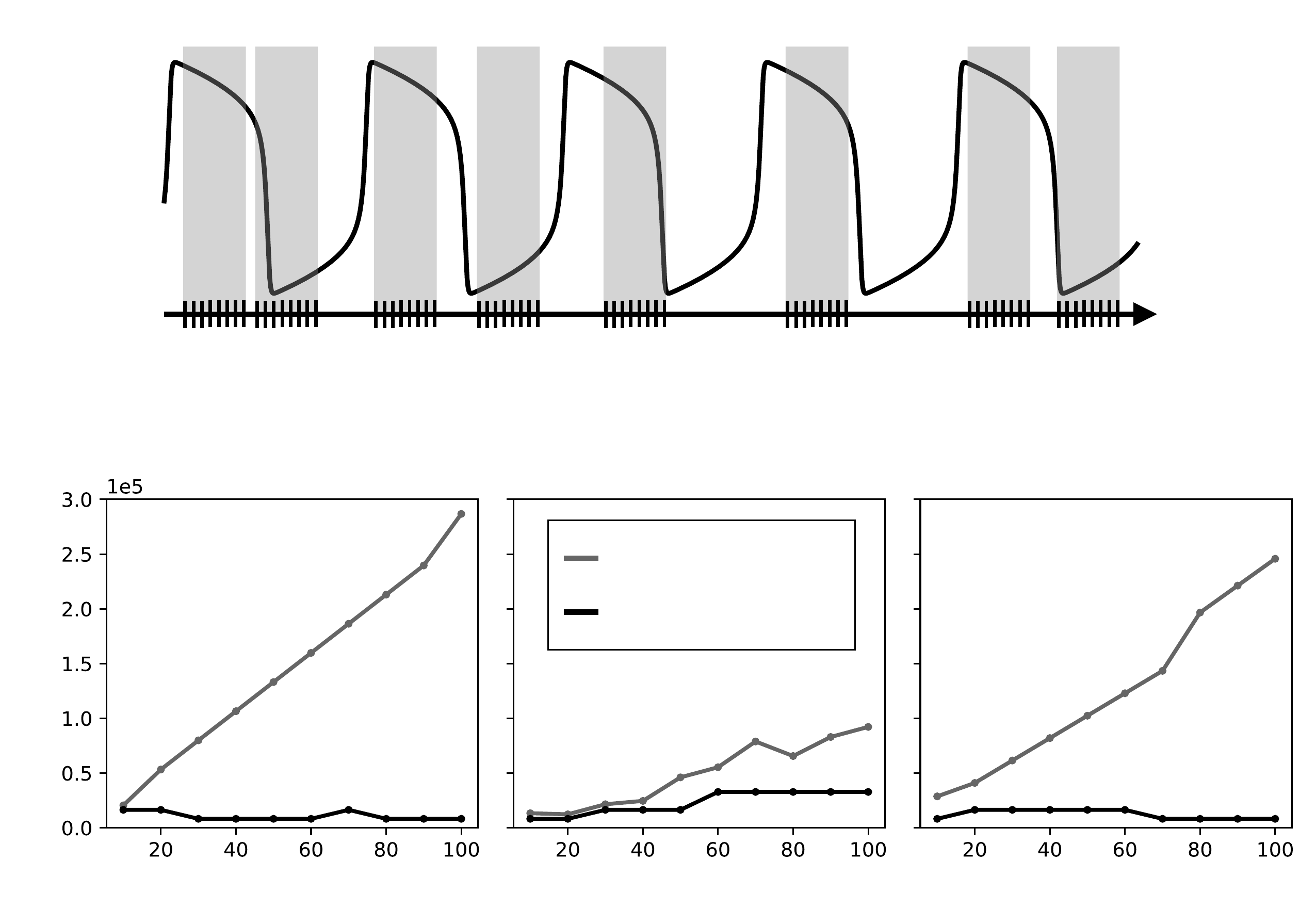}
  \put(3,68){\color{header1}\HLtext\normalsize (a) Burst sampling}
  \put(3,39.5){\color{header1}\HLtext\normalsize (b) Sampling requirements}
  \put(12,34){\small \parbox{8em}{\centering Fast and Slow Van der Pol}}
  \put(45,34){\small \parbox{8em}{\centering Slow Van der Pol, Fast Lorenz}}
  \put(75,34){\small \parbox{8em}{\centering Fast Van der Pol, Slow Lorenz}}
  \put(47,26.5){\small uniform sampling}
  \put(47,22){\small burst sampling}
  \put(1,7){\small \rotatebox{90}{\parbox{10em}{\centering samples required}}}
  \put(47,0){\small frequency ratio}
  \end{overpic}
  \caption{Schematic and results of our burst sampling method for applying SINDy to multiscale systems. (a) Illustration of burst sampling. Samples are collected in short bursts at a fine resolution, spaced out randomly over a long duration. This allows us to sample for a long time with a reduced effective sampling rate (as compared with naive uniform sampling) without increasing the time step. (b) Comparison of uniform sampling and burst sampling on systems with multiple time scales. For each method we look at how many samples are required for SINDy to identify the correct model at different frequency ratios between the fast and slow time scales. With uniform sampling, the total sampling requirement increases significantly as the frequency ratio between the fast and slow systems increases. Burst sampling significantly reduces this scaling effect, allowing SINDy to scale efficiently to multiscale systems.}
  \label{fig:burst_sampling}
\end{figure}
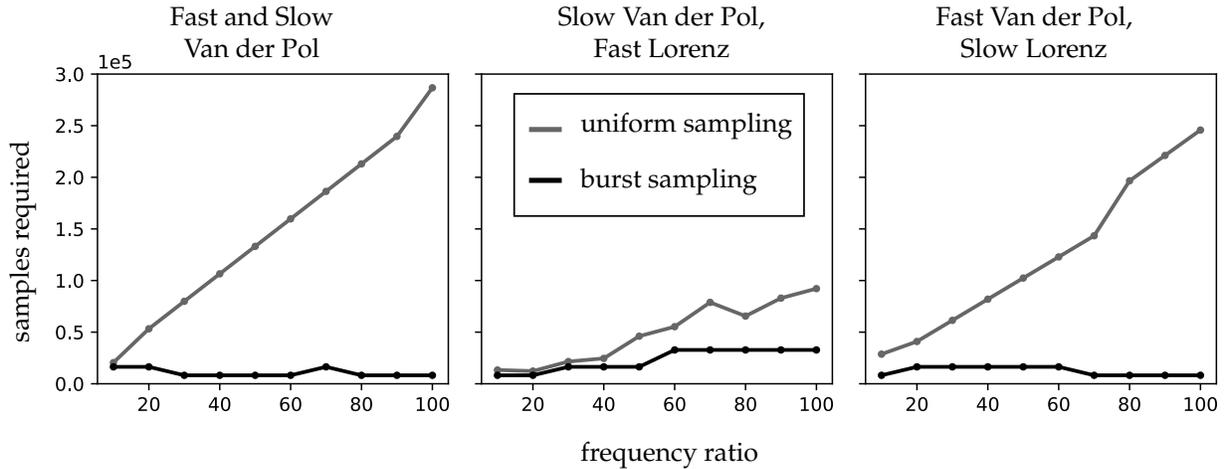

We have established that SINDy can identify uniscale dynamical systems with relatively little data. However, many systems of interest contain coupled dynamics at multiple time scales. Identifying the dynamics across scales would require sampling at a sufficiently fast rate to capture the fast dynamics while also sampling for a sufficient duration to observe the slow dynamics. Assuming we collect samples at a uniform rate, this leads to an increase in the amount of data required as the time scales of the coupled system separate. In this section, we introduce a sampling method that overcomes this issue, allowing SINDy to scale efficiently to multiscale problems.

We consider coupled systems with two distinct time scales. In particular, we look at nonlinear systems with linear coupling:
\begin{subequations}
\begin{align}
\tau_{\text{fast}} \dot{\mathbf{u}} &= \mathbf{f}(\mathbf{u}) + \mathbf{C} \mathbf{v} \nonumber \\
\tau_{\text{slow}} \dot{\mathbf{v}} &= \mathbf{g}(\mathbf{v}) + \mathbf{D} \mathbf{u}. \nonumber
\end{align}
\end{subequations}
In this system $\mathbf{u}(t)\in\mathbb{R}^n$ is the set of variables that comprise the fast dynamics, and $\mathbf{v}(t)\in\mathbb{R}^l$ represents the slow dynamics. The linear coupling is determined by our choice of $\mathbf{C}\in\mathbb{R}^{n\times l},\mathbf{D}\in\mathbb{R}^{l\times n}$, and time constants $\tau_\text{fast},\tau_\text{slow}$, which determine the frequency of the fast and slow dynamics. We consider three example systems: two coupled Van der Pol oscillators, a slow Van der Pol oscillator coupled with a fast Lorenz system, and a fast Van der Pol oscillator coupled with a slow Lorenz system. In order to understand the effect of time scale separation, we consider the frequency ratio $F=T_\text{slow}/T_\text{fast}$ between the coupled systems; $T_\text{slow},T_\text{fast}$ are the approximate ``periods'' of the slow and fast systems, defined for Lorenz and Van der Pol in Sections~\ref{sec:lorenz}, \ref{sec:vdp} respectively. We assess how much data is required to discover the system as $F$ increases.

For each of the three example multiscale systems, we assess the performance of SINDy using the same process outlined in Section~\ref{sec:uniscale1}. Using naive uniform sampling, the data requirement increases approximately linearly with the frequency ratio $F$ for all three systems (see Figure~\ref{fig:burst_sampling}). This means that the sampling requirement is extremely high for systems where the frequency scales are highly separated (large $F$). Because the tasks of data collection and of fitting models to large data sets can be computationally expensive, reducing the data required by SINDy is advantageous.

We introduce a sampling strategy to address this issue, which we refer to as burst sampling. By using burst sampling, the data requirement for SINDy stays approximately constant as time scales separate ($F$ increases). The top panel of Figure~\ref{fig:burst_sampling} shows an illustration of our burst sampling method. The idea is to maintain a small step size but collect samples in short bursts spread out over a long duration. This reduces the effective sampling rate, which is particularly useful when we are limited by bandwidth in the number of samples we can collect. By maintaining a small step size, we still observe the fine detail of the dynamics, and we can get a more accurate estimate of the derivative. However by spacing out the samples in bursts, our data also captures more of the slow dynamics than would be observed with the same amount of data collected at a uniform sampling rate. We thus reduce the effective sampling rate without degrading our estimates of the derivative or averaging out the fast time scale dynamics.

Employing burst sampling requires the specification of a few parameters. We must select a step size $\Delta t$, burst size (number of samples in a burst), duration over which to sample, and total number of bursts to collect. We fix a single step size and consider the effect of burst size, duration, and number of bursts. In general, we find that smaller burst sizes give better performance. We also find that while it is important to have a sufficiently long duration (on the order of 1-2 periods of the slow dynamics), increasing duration beyond this does not have a significant effect on performance. Therefore, for the rest of our analysis we fix a burst size of $8$ and a sampling duration of $2T_\text{slow}$. We then adjust the total number of bursts collected in order to control the effective sampling rate.

Another important consideration is how to space out the bursts of samples. In particular we must consider the potential effects of aliasing, which could reduce the effectiveness of this method if bursts are spaced to cover the same portions of the attractor. To address this issue, we introduce randomness into our decision of where to place bursts. In streaming data applications, this can be handled in a straightforward manner by selecting burst collection times as Poisson arrival times, with the rate of the Poisson process chosen so that the expected number of samples matches our desired overall sampling rate. For the purpose of testing our burst sampling method, we do something slightly different. In order to understand the performance of burst sampling, we need to perform repeated trials with different choices of burst locations. We also need to ensure that for each choice of parameters, the training data for each trial has a constant number of samples and covers the same duration. If burst locations are selected completely randomly, it is likely that some will overlap which would reduce our number of overall samples. To address this, we start with evenly spaced bursts and sample offsets for each, selected from a uniform distribution that allows each burst to shift left or right while ensuring that none of the bursts overlap. In this manner we obtain the desired number of samples at each trial, but we introduce enough randomness into the selection of burst locations to account for aliasing.

In Figure~\ref{fig:burst_sampling} we show the results of burst sampling for our three example systems. For all three systems, the number of samples required by SINDy remains approximately constant as the frequency ratio $F$ increases. To determine the number of samples required at a given frequency ratio, we fix an effective sampling rate and run repeated trials with burst locations selected randomly as described above. For each trial, we apply SINDy to the sampled data and determine if the correct nonzero coefficients were identified, as in Section~\ref{sec:uniscale1}. We run 100 trials at each effective sampling rate and find the minimum rate for which the system was correctly identified in all trials; this determines the required effective sampling rate. In Figure~\ref{fig:burst_sampling}, we plot the total number of samples collected with this effective sampling rate.

\section{Incomplete measurements of the state space and latent variables}
\label{sec:incomplete_measurements}

The results so far have assessed the performance of the SINDy algorithm on coupled multiscale systems with limited full-state measurements in time. One limitation of SINDy is that it requires knowledge of the full underlying state space variables that govern the behavior of the system of interest. In many real-world applications, some governing variables may be completely unobserved or combined into mixed observations. With multiscale systems in particular, a single observation variable may contain dynamics from multiple time scales. We therefore need methods for understanding multiscale dynamics that do not require full state measurements.

One typical approximation technique for modeling dynamical systems relies on attempts to linearize the nonlinear dynamics. Linear models for dynamics have many advantages, particularly for control and prediction.  Koopman analysis is an emerging data-driven modeling tool for dynamical systems.  First proposed in 1931 \cite{koopman_hamiltonian_1931}, it has experienced a recent resurgence in interest and development \cite{mezic_spectral_2005,budisic_applied_2012,mezic_analysis_2013}.   
Here we assume that we are working with a discrete-time version of the dynamical system in (\ref{eq:dynamical_system}):
\begin{equation}
  \mathbf{x}_{k+1} = \mathbf{F}(\mathbf{x}_k) = \mathbf{x}_k + \int_{k\Delta t}^{(k+1)\Delta t} \mathbf{f}(\mathbf{x}(\tau))\mathrm{d}\tau.
\end{equation}
The Koopman operator $\mathcal{K}$ is an infinite-dimensional linear operator on the Hilbert space of measurement functions of the states $\mathbf{x}_k$. Given measurements $g(\mathbf{x}_k)$, the Koopman operator is defined by
\begin{equation}
  \mathcal{K}g \triangleq g \circ \mathbf{F} \quad \Rightarrow \quad \mathcal{K}g(\mathbf{x}_k) = g(\mathbf{x}_{k+1}).
\end{equation}
Thus the Koopman operator maps the system forward in time in the space of measurements.

While having a linear representation of the dynamical system is advantageous, the Koopman operator is infinite-dimensional and obtaining finite-dimensional approximations is difficult in practice. In order to obtain a good model, we seek a set of measurements that form a Koopman invariant subspace~\cite{brunton_koopman_2016}. Dynamic mode decomposition (DMD)\cite{schmid_dynamic_2010} is one well-known method for approximating the Koopman operator~\cite{Rowley2009jfm,Tu2014jcd,kutz_dynamic_2016}. DMD constructs a linear mapping satisfying
\begin{equation}
  \mathbf{x}_{k+1} \approx \mathbf{A}\mathbf{x}_k. \nonumber
\end{equation}
However as might be expected, DMD does not perform well for strongly nonlinear systems. Extended DMD and kernel DMD are two methods that seek to resolve this issue by constructing a library of nonlinear measurements of $\mathbf{x}_k$ and finding a linear operator that works on these measurements \cite{williams_datadriven_2015,williams_kernel-based_2015}. However these methods can be computationally expensive, and it is not guaranteed that the selected measurements will form a Koopman invariant subspace~\cite{brunton_koopman_2016} unless results are rigorously cross-validated, as in the equivalent variational approach of conformation dynamics (VAC) approach~\cite{noe2013variational,nuske2014jctc}.  
Alternatively, one can find judicious choices of the nonlinear measurements that transform the underlying dynamics from a strongly nonlinear system to a weakly nonlinear system~\cite{kutz2016koopman}.
A review of the DMD algorithm and its many applications can be found in Ref.~\cite{kutz_dynamic_2016}

The recent Hankel alternative view of Koopman (HAVOK) method constructs an approximation to the Koopman operator by relying on the relationship between the Koopman operator and the Takens embedding \cite{takens_detecting_1981,brunton_chaos_2017}; delay coordinates were previously used to augment the rank in DMD~\cite{Tu2014jcd,brunton2016extracting} and the connection to Koopman theory has been strengthened~\cite{Arbabi2016arxiv,Das2017arxiv} following the original HAVOK paper~\cite{brunton_chaos_2017}. Measurements of the system are formed into a Hankel matrix, which is created by stacking delayed measurements of the system:
\begin{equation}
  \mathbf{H} = \left( \begin{array}{cccc}
    \mathbf{x}_1 & \mathbf{x}_2 & \cdots & \mathbf{x}_p \\
    \mathbf{x}_2 & \mathbf{x}_3 & \cdots & \mathbf{x}_{p+1} \\
    \vdots & \vdots & \ddots & \vdots \\
    \mathbf{x}_q & \mathbf{x}_{q+1} & \cdots & \mathbf{x}_m
  \end{array} \right).
\end{equation}
A number of other algorithms make use of the Hankel matrix, including the eigensystem realization algorithm (ERA) \cite{juang_eigensystem_1985,broomhead_time-series_1989,le_clainche_higher_2017}. By taking a singular value decomposition (SVD) of the Hankel matrix, we are able to obtain dominant time-delay coordinates that are approximately invariant to the Koopman operator \cite{brunton_chaos_2017}.  Thus the time-delay embedding provides a new coordinate system in which the dynamics are linearized.

In \cite{brunton_chaos_2017} the focus is on chaotic systems, and an additional nonlinear forcing term is included in the HAVOK model to account for chaotic switching or bursting phenomena. Here we focus on quasiperiodic systems and show that the linear model found by HAVOK is sufficient for reconstruction and long-term prediction, with no need for a forcing term. HAVOK has the advantage that the discovered models are linear and require no prior knowledge of the true governing variables. By eliminating the need for the nonlinear forcing term, we obtain deterministic, closed-form models. In this section we assess the performance of HAVOK on uniscale dynamical systems and introduce two strategies for scaling the method to problems with multiple time scales.

\subsection{Uniscale dynamical systems}
\label{sec:uniscale_incomplete}

\begin{figure}[t]
  \centering
  \begin{overpic}[width=\linewidth]{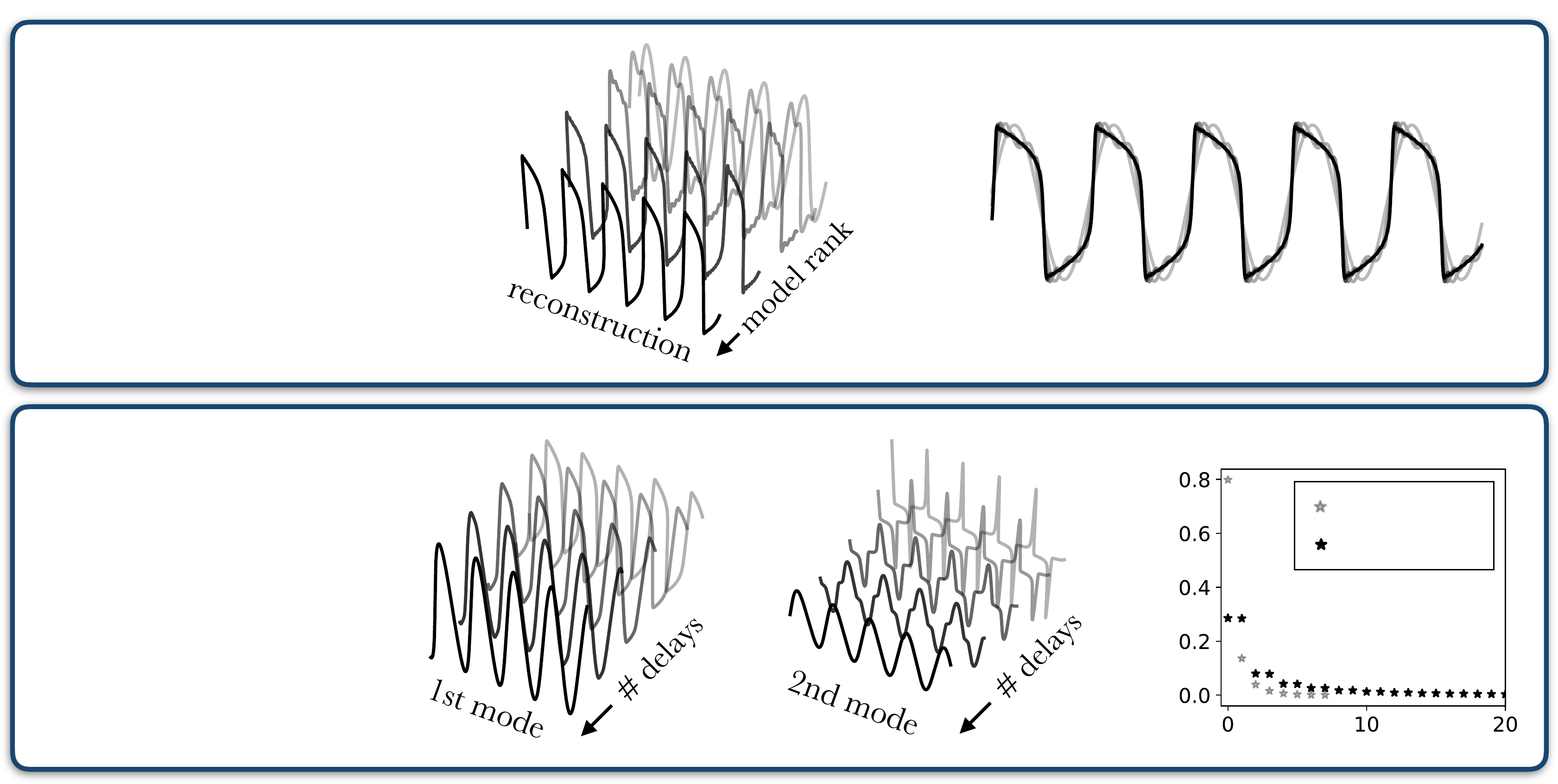}
  \put(3,37){\parbox{8em}{\centering \color{header1}\HLtext Effect of Rank on HAVOK Model}}
  \put(3,12){\parbox{8em}{\centering \color{header1}\HLtext Effect of Delays on HAVOK Model}}
  \put(80.5,21.5){\footnotesize Singular Values}
  \put(87,17.5){\scriptsize 8 delays}
  \put(87,14.8){\scriptsize 64 delays}
  \put(87.5,1.5){\scriptsize $k$}
  \put(72.5,7){\scriptsize \rotatebox{90}{\parbox{5em}{\centering $\sigma_k/\sum_i \sigma_i$}}}
  \end{overpic}
  \vspace{-.2in}
  \caption{Effects of rank and delays on HAVOK model of a Van der Pol oscillator. In the top panel, we show how the choice of model rank affects the HAVOK model. As the rank of the model increases, the model goes from being a sinusoidal linear oscillator to closely matching the Van der Pol dynamics. In the bottom panel, we show how the number and duration of the delay coordinates affects the HAVOK modes and singular values. HAVOK models are linear, and linear dynamics consist of the superposition of sinusoids; thus to obtain a good model the modes found by taking an SVD of $\mathbf{H}$ should resemble Fourier modes (in time). With few delays the modes are highly nonlinear, meaning that our (linear) HAVOK model cannot accurately capture the dynamics. As the number of delays increases the modes appear more like Fourier modes, allowing us to construct a good linear model. The right panel shows the singular values for two models, one with 8 delays and one with 64 delays (only the first 20 singular values are shown). Using 8 delays admits a low-rank model with nonlinear modes; most of the energy is captured in the first mode. Using 64 delays admits linear modes, and the energy is more evenly shared among the first several modes; therefore several modes must be included to get an accurate model.}
  \label{fig:5a}
\end{figure}

\begin{figure}[t]
  \centering
  \begin{overpic}[width=\linewidth]{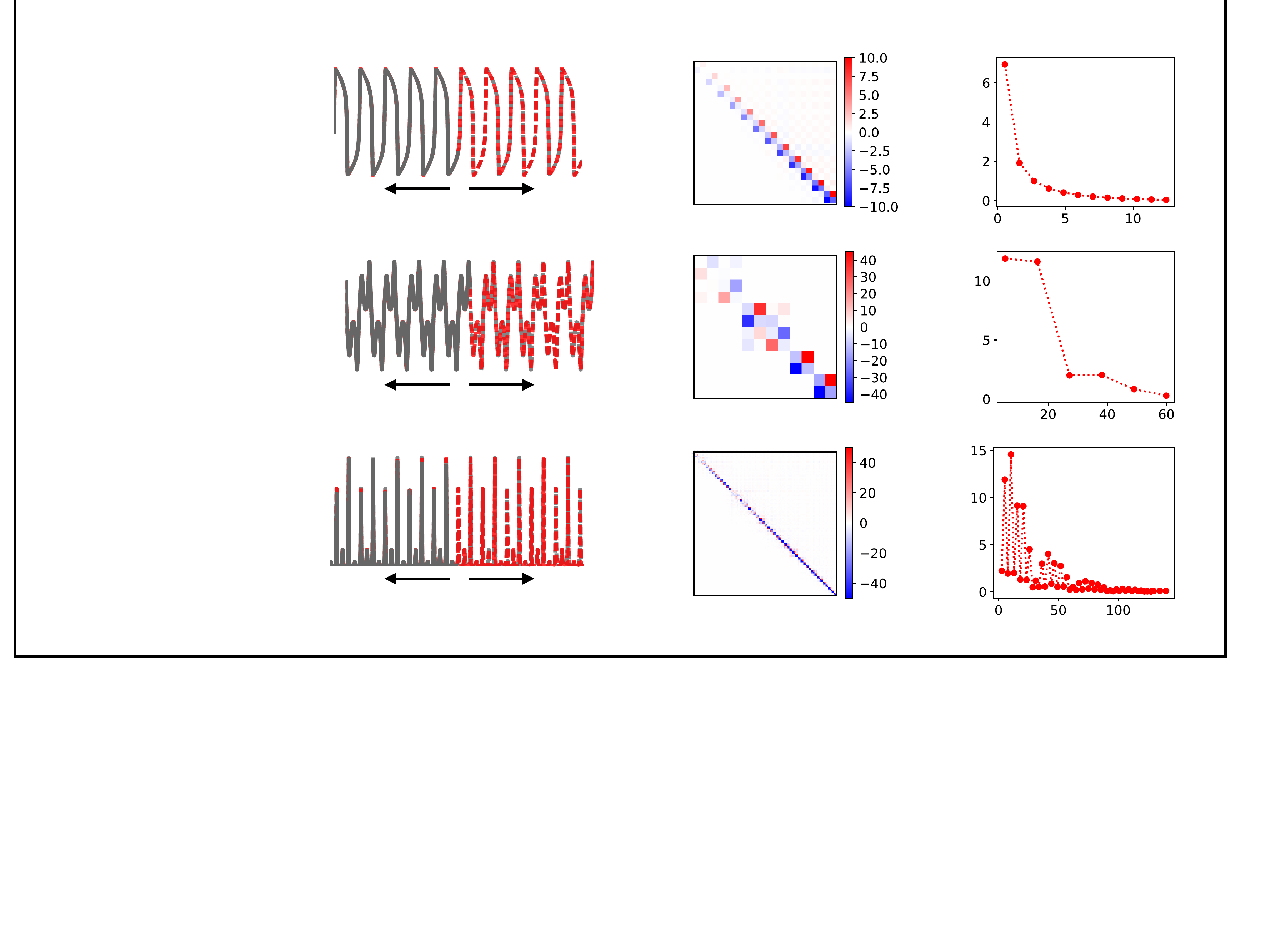}
  \put(7,53){\color{header1}\HLtext System}
  \put(27.5,53){\color{header1}\HLtext HAVOK Model}
  \put(52,53){\color{header1}\HLtext Dynamics Matrix}
  \put(79,53){\color{header1}\HLtext Power Spectrum}
  \put(0,47){\parbox{10em}{\centering Van der Pol (\ref{eq:vdp})}}
  \put(0,30){\parbox{10em}{\centering Lorenz (periodic) (\ref{eq:lorenz})}}
  \put(0,13){\parbox{10em}{\centering Rossler (periodic) (\ref{eq:rossler})}}
  \put(0,42){\parbox{10em}{\centering\footnotesize $\mu=5$ \\ variable modeled: $x$ \\ model rank: 24}}
  \put(0,25){\parbox{10em}{\centering\footnotesize $\sigma=10,\ \beta=\frac{8}{3},\ \rho=160$ \\ variable modeled: $x$ \\ model rank: 12}}
  \put(0,8){\parbox{10em}{\centering\footnotesize $a=0.1,\ b=0.1,\ c=8.5$ \\ variable modeled: $z$ \\ model rank: 105}}
  \put(29,36.5){\scriptsize training}
  \put(38,36.5){\scriptsize test}
  \put(29,20){\scriptsize training}
  \put(38,20){\scriptsize test}
  \put(29,4){\scriptsize training}
  \put(38,4){\scriptsize test}
  \put(84.5,1){\scriptsize frequency}
  \end{overpic}
  \vspace{-.3in}
  \caption{Hankel alternative view of Koopman (HAVOK) models for three example systems. We train HAVOK models for a single variable of each system. For each example we show the model reconstruction of the dynamics, the matrix defining the linear model of the dynamics, and the power spectrum of the DMD modes.}
  \label{fig:5b}
\end{figure}
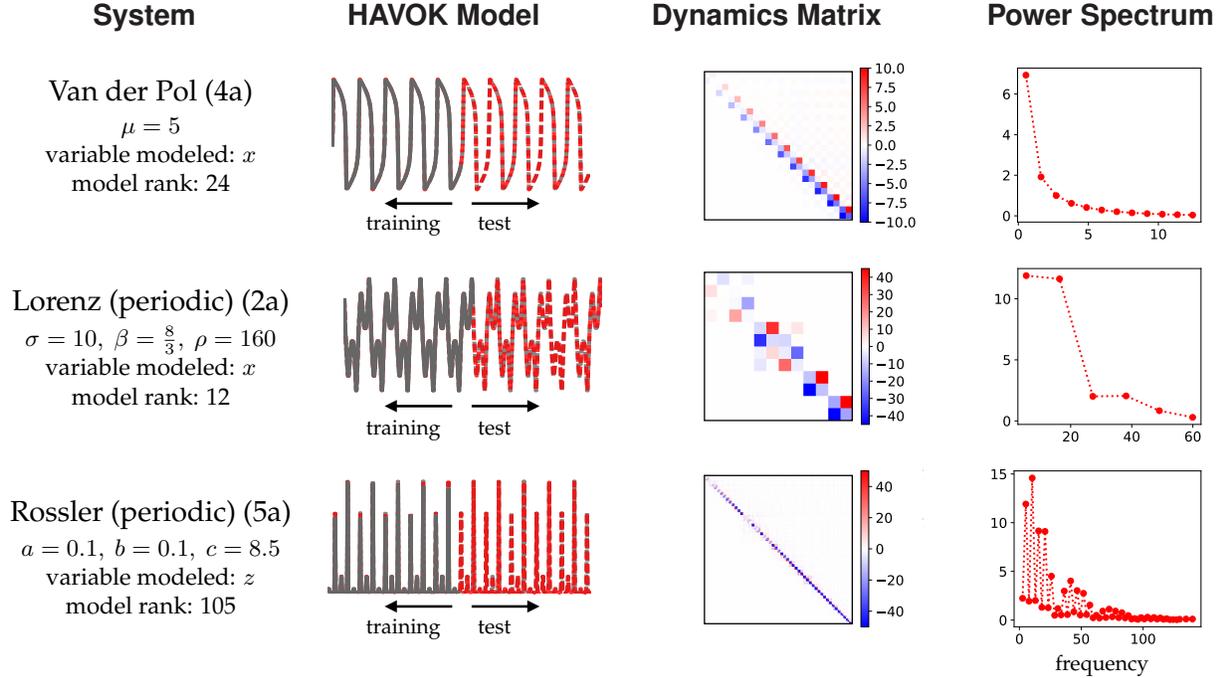

We start by focusing on quasiperiodic systems with a single time scale. To apply the HAVOK method, we form two shift-stacked matrices that are analogous to the snapshot matrices typically formed in DMD~\cite{kutz_dynamic_2016}:
\begin{equation}
    \mathbf{H} = \left( \begin{array}{cccc}
    \mathbf{x}_1 & \mathbf{x}_2 & \cdots & \mathbf{x}_{m-q} \\
    \mathbf{x}_2 & \mathbf{x}_3 & \cdots & \mathbf{x}_{m-q+1} \\
    \vdots & \vdots & \ddots & \vdots \\
    \mathbf{x}_q & \mathbf{x}_{q+1} & \cdots & \mathbf{x}_{m-1}
  \end{array} \right), \quad
  \mathbf{H}' = \left( \begin{array}{cccc}
    \mathbf{x}_2 & \mathbf{x}_3 & \cdots & \mathbf{x}_{m-q+1} \\
    \mathbf{x}_3 & \mathbf{x}_4 & \cdots & \mathbf{x}_{m-q+2} \\
    \vdots & \vdots & \ddots & \vdots \\
    \mathbf{x}_{q+1} & \mathbf{x}_{q+2} & \cdots & \mathbf{x}_m
  \end{array} \right).
  \label{eq:havok_snapshots}
\end{equation}
We then perform the standard DMD algorithm on the shift-stacked matrices. A rank-$r$ DMD model gives us a set of $r$ eigenvalues $\lambda_k$, modes $\bm{\phi}_k \in \mathbb{C}^n$, and amplitudes $b_k$. By rewriting the discrete-time eigenvalues $\lambda_k$ as $\omega_k = \ln(\lambda_k)/\Delta t$, we obtain a linear model
\begin{equation}
  \mathbf{x}(t) = \sum_{k=1}^r \bm{\phi}_k\exp(\omega_k t) b_k.
  \label{eq:dmd_model}
\end{equation}
Note that this linear model provides a prediction of the behavior at any time $t$, without requiring the system to be fully simulated up to that time.

In Figure~\ref{fig:5b} we show HAVOK models for the periodic Van der Pol, Rossler, and Lorenz systems. In each example we take observations of a single variable of the system and time delay embed to build a linear HAVOK model. At the selected parameter values, each of these systems exhibits regular periodic behavior. We also show the power spectrum of each HAVOK model by plotting the model frequencies $\mathrm{Im}(\omega_k)$ and their respective amplitudes $|b_k|$. 

When applying HAVOK to data, we must consider the sampling rate, number of delays $q$, and rank of the model $r$. In order to obtain a good model, it is imperative that (1) there are enough delay coordinates to provide a model of sufficient rank and (2) the delay duration $D=(q-1)\Delta t$ must be large enough to capture a sufficient duration of the oscillation. To obtain a good model, we need the modes found by taking an SVD of $\mathbf{H}$ to resemble Fourier modes. This is due to the fact that HAVOK builds a linear dynamical model: linear dynamics without damping consists of the superposition of pure tone sinusoids in time. Highly nonlinear modes will not be well captured by a linear HAVOK model. In Figure~\ref{fig:5a} we show how the modes for a Van der Pol oscillator go from nonlinear to Fourier modes as we increase the number of time delays (which simultaneously increases the delay duration $D$). As a rule of thumb, for a system with period $T$ we choose our delays such that $D=T$. Note that with a high sampling rate, this could necessitate using a significant number of delay coordinates, making the snapshot matrices impractically large. Two options for dealing with this are (1) downsampling the data and (2) spacing out the delays; this is discussed further in Section~\ref{sec:multiscale2_spacing}. For all examples in Figure~\ref{fig:5b}, we use $q=128$ delays and take $\Delta t=T/(q-1)$. We fit the model using training data from a single trajectory covering $5$ periods of oscillation ($m\Delta t=5T$).

 Choosing the proper rank of the HAVOK model is another important consideration. Increasing the rank of the model adds additional frequencies, which in general can lead to a more accurate reconstruction. The top panel of Figure~\ref{fig:5a} shows how increasing rank affects a HAVOK model of a Van der Pol oscillator. A rank-2 model looks like a linear oscillator with the same dominant frequency as the Van der Pol oscillator. As we increase the rank, the model becomes more like a Van der Pol oscillator, but we observe Gibbs phenomenon near the sharp transitions in the dynamics. With a rank-32 model we get a very close reconstruction of the Van der Pol behavior. In order to obtain a model of rank $r$, we must choose the number of delays $q$ such that $q>r$. Assuming the delay embedding is constructed adequately for capturing the dynamics ($m,D$ sufficiently large and $\Delta t$ sufficiently small), the rank of the model could be selected using traditional Pareto front analysis.

\begin{figure}[t]
  \centering
  \begin{overpic}[width=0.9\linewidth]{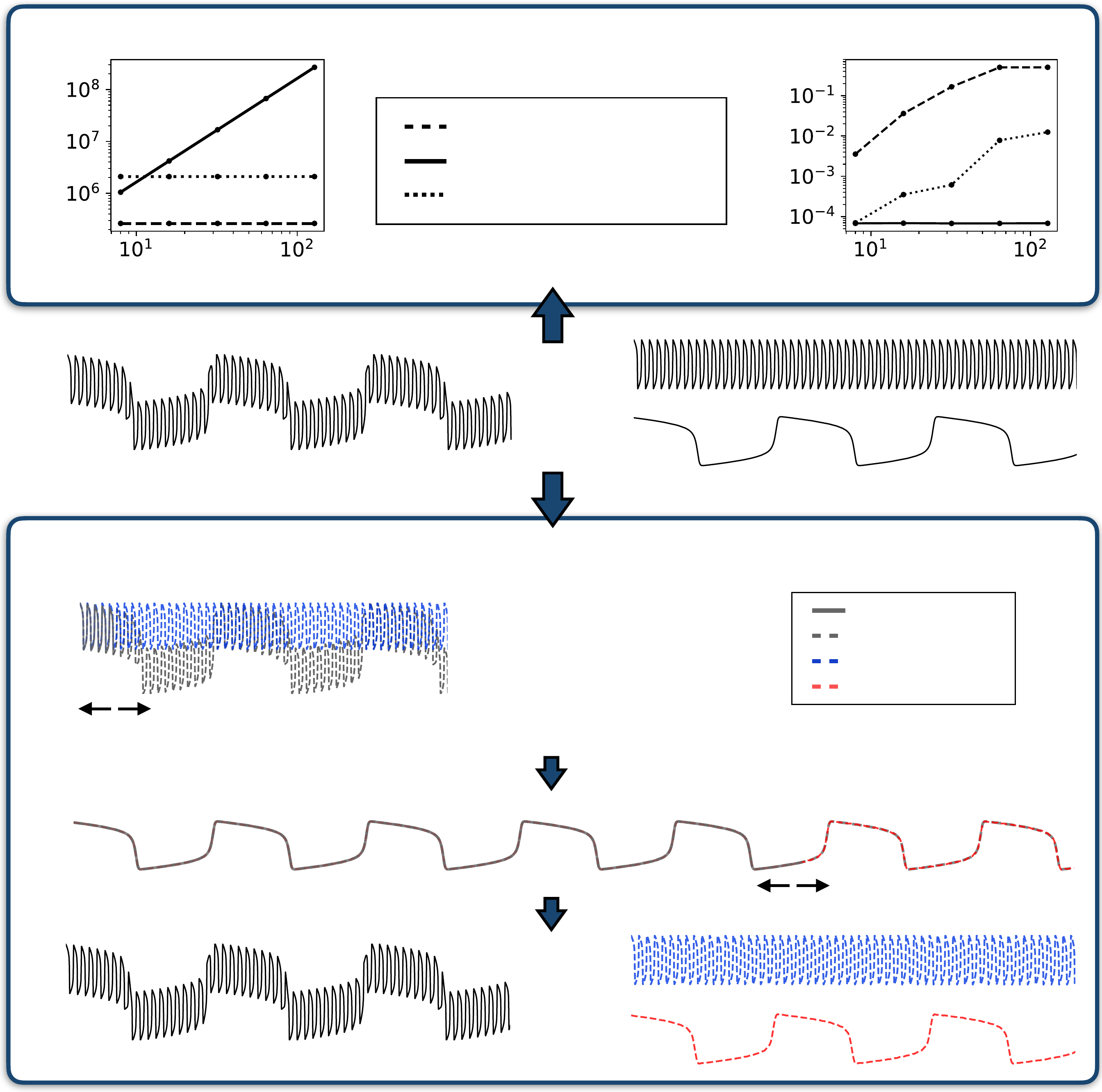}
  \put(38,95.5){\color{header1}\HLtext (a) Delay Spacing}
  \put(12.5,73){\footnotesize frequency ratio}
  \put(3.2,77){\small \rotatebox{90}{\parbox{7em}{\centering $\texttt{numel}(\mathbf{H})$}}}
  \put(68.5,77){\small \rotatebox{90}{\parbox{7em}{\centering RMSE}}}
  \put(79,73){\footnotesize frequency ratio}
  \put(42,87.5){\footnotesize no spacing, large $\Delta t$}
  \put(42,84){\footnotesize no spacing, small $\Delta t$}
  \put(42,81){\footnotesize spacing, small $\Delta t$}
  \put(38,49){\color{header1}\HLtext (b) Iterative Modeling}
  \put(1.5,62){\small $y(\mathbf{x})$}
  \put(48,62){\small $=$}
  \put(51,65){\small $x_{\text{fast}}$}
  \put(52,62){\small $+$}
  \put(51,58.5){\small $x_{\text{slow}}$}
  \put(3.5,32.5){\scriptsize training}
  \put(11,32.5){\scriptsize test}
  \put(65,16){\scriptsize training}
  \put(72.5,16){\scriptsize test}
  \put(2,47){\color{header1}\HLtext \footnotesize Model fast dynamics}
  \put(2,45){\scriptsize high sampling rate, short duration}
  \put(2,29){\color{header1}\HLtext \footnotesize Pull out fast dynamics, model slow dynamics}
  \put(2,27){\scriptsize low sampling rate, long duration}
  \put(2,15){\color{header1}\HLtext \footnotesize Combine fast and slow models}
  \put(78.5,43.5){\scriptsize training data}
  \put(78.5,41){\scriptsize held out data}
  \put(78.5,38.5){\scriptsize fast model}
  \put(78.5,36){\scriptsize slow model}
  \put(2,7){\small $y(\mathbf{x})$}
  \put(47,7){\small $\approx$}
  \put(49.5,7){\parbox{3em}{\centering\footnotesize fast model \\ + \\ slow model}}
  \end{overpic}
  \vspace*{-.1in}
  \caption{Strategies for applying HAVOK to multiscale systems, using an example system of a fast and a slow Van der Pol oscillator. (a) Spacing out the rows and columns of the HAVOK snapshot matrices while maintaining a small time step allows us to make a flexible trade-off between computational cost and model error. We show how the choice of time step, number of delay coordinates, and spacing affects the size of the embedding matrices (left) and model error (right). (b) Schematic of our iterative method for applying HAVOK to multiscale systems. First, HAVOK is applied to training data from a very short recording period to capture just the fast dynamics. Next, the data is sampled over a longer duration at a much lower rate. The known fast dynamics are subtracted out and the resulting data is used to create a HAVOK model for the slow dynamics. The fast and slow models can be combined to model the full system.}
  \label{fig:5}
\end{figure}
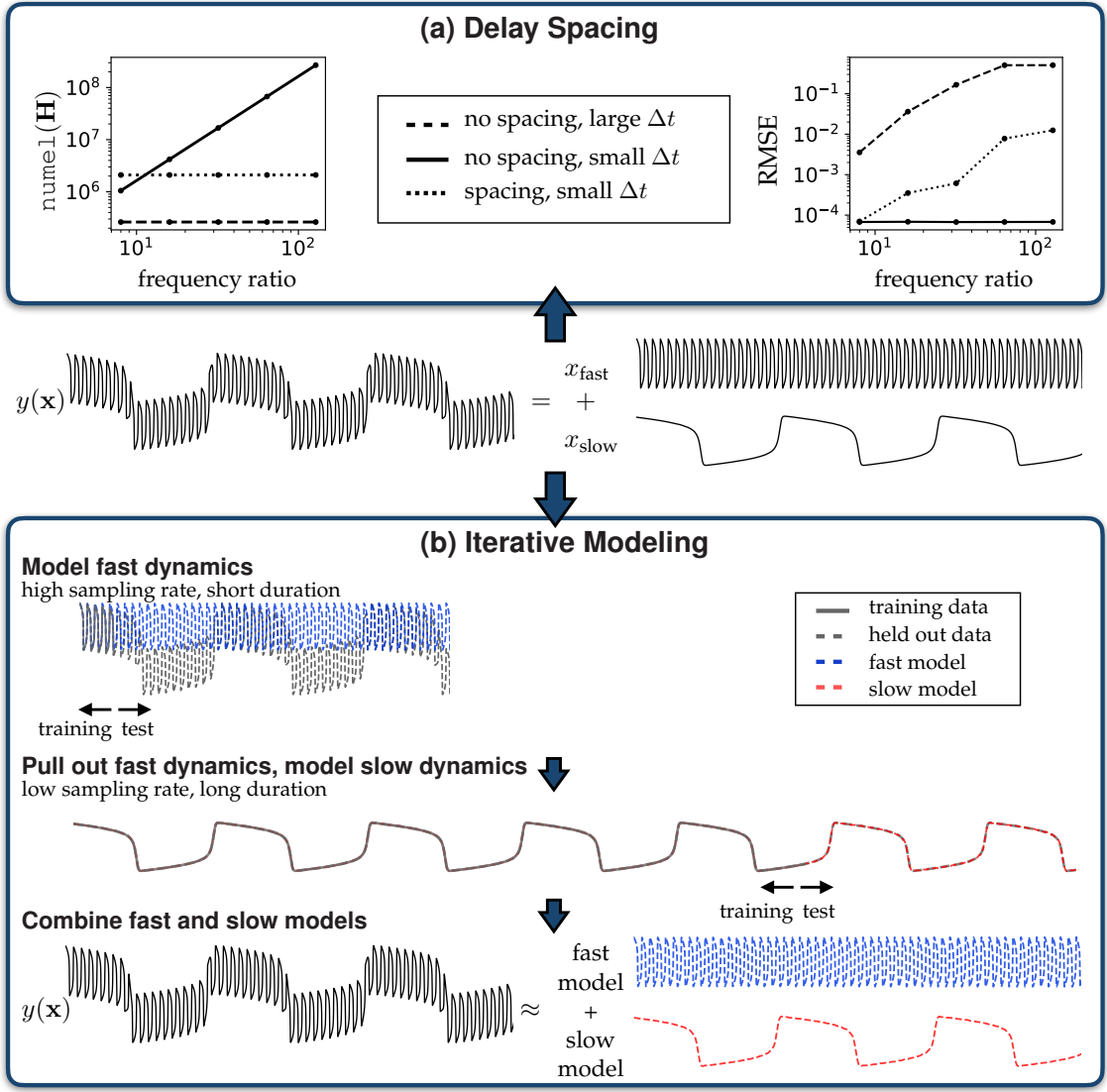

Our HAVOK model is constructed by taking an SVD of the delay embedding matrix, which consists of samples from a finite sampling period. This introduces some bias into the calculation of the eigenvalues of the systems. In particular, for periodic and quasiperiodic systems where we are trying to capture the behavior on the attractor, we do not want our model to have blowup or decay. Thus the continuous eigenvalues $\omega_k$ should have zero real part. In general, when we apply HAVOK we obtain eigenvalues with a small but nonzero real part. To deal with this issue, we simply set the real part of all eigenvalues to zero after fitting the HAVOK model. This approach is also taken in \cite{le_clainche_higher_2017}. Forcing all eigenvalues to be purely imaginary allows the model to predict long term behavior of the system without blowup or decay. While we find this sufficient to produce good models for the example systems studied, an alternative solution would be to use a method such as optimized DMD with an imposed constraint on the eigenvalues \cite{askham_variable_2018}.

\subsection{Multiscale dynamical systems}
\label{sec:multiscale2}

It is still feasible to apply HAVOK as described above to many systems with multiple time scales. As with SINDy, the time scale separation requires that the amount of data acquired by the algorithm increase to account for both the fast and slow dynamics so that the problem becomes more computationally expensive. The use of delay coordinates compounds this issue beyond just a blowup in the number of samples required: not only do we need a sufficiently small time step to capture the fast dynamics, we also need a sufficient number of delays $q$ so that our delay duration $D=(q-1)\Delta t$ covers the dynamics of the slow time scale. Thus the size of our snapshot matrices scales as $F^2$ with the frequency ratio $F$ between fast and slow time scales. For large frequency separations and systems with high-dimensional state variables, this results in an extremely large delay embedding matrix that could make the problem computationally intractable. In this section, we discuss two strategies for applying HAVOK to problems with multiple time scales.

\subsubsection{Method 1: Delay spacing}
\label{sec:multiscale2_spacing}

Our first approach for applying HAVOK to multiscale problems is a simple modification to the form of the delay embedding matrices. Rather than using the standard embedding matrices given in \eqref{eq:havok_snapshots}, we introduce row and column spacings $d,c$ and create the following snapshot matrices:

\begin{align}
    \mathbf{H} &= \left( \begin{array}{cccc}
    \mathbf{x}_1 & \mathbf{x}_{c+1} & \cdots & \mathbf{x}_{(p-1)c + 1} \\
    \mathbf{x}_{d+1} & \mathbf{x}_{d+c+1} & \cdots & \mathbf{x}_{d + (p-1)c + 1} \\
    \vdots & \vdots & \ddots & \vdots \\
    \mathbf{x}_{(q-1)d+1} & \mathbf{x}_{(q-1)d+c+1} & \cdots & \mathbf{x}_{(q-1)d + (p-1)c + 1}
  \end{array} \right), \quad \nonumber \\
  \mathbf{H}' &= \left( \begin{array}{cccc}
    \mathbf{x}_2 & \mathbf{x}_{c+2} & \cdots & \mathbf{x}_{(p-1)c + 2} \\
    \mathbf{x}_{d+2} & \mathbf{x}_{d+c+2} & \cdots & \mathbf{x}_{d + (p-1)c + 2} \\
    \vdots & \vdots & \ddots & \vdots \\
    \mathbf{x}_{(q-1)d+2} & \mathbf{x}_{(q-1)d+c+2} & \cdots & \mathbf{x}_{(q-1)d + (p-1)c + 2}
  \end{array} \right). \nonumber
  \label{eq:havok_snapshots_spacing}
\end{align}
The number of columns $p$ is determined by the total number of samples. This form of delay embedding matrix separates out both the rows (delay coordinates) and columns (samples) without increasing the time step used in creating the HAVOK model. Each shift-stacked copy of the data advances in time by $d\Delta t$, which means the delay duration $D$ is now multiplied by a factor of $d$, giving $D=(q-1)d\Delta t$. If we keep $q$ fixed and increase $d$, our delay embedding captures more of the dynamics without increasing the number of rows in the embedding matrices. Similarly, we construct $\mathbf{H}$ and $\mathbf{H}'$ so that each column advances in time by $c\Delta t$, meaning we reduce the number of columns in the snapshot matrices by a factor of $c$ without reducing the total sampling duration. The key point is that we keep the time step small by constructing $\mathbf{H}'$ so that it only advances each data point by time step $\Delta t$ from the snapshots in $\mathbf{H}$. This allows the HAVOK model to account for the fast time scale.

We consider multiscale systems with two time scales and assess performance as the frequency ratio between the time scales grows. In particular, as the frequency ratio increases we are concerned with both the size of the snapshot matrices and the prediction error of the respective model. We compare three different models: a standard HAVOK model with a large time step, a standard HAVOK model with a small time step, and a HAVOK model built with row and column spacing in the snapshot matrices. As discussed in Section~\ref{sec:uniscale_incomplete}, the most important considerations in building a HAVOK model are selecting the model rank $r$ and using a sufficiently large delay duration $D$. In our comparison we keep $r=100$ fixed and set $D$ equal to the period of the slow dynamics.

As an example we consider a system of two Van der Pol oscillators with parameter $\mu=5$ and periods $T_\text{fast},T_\text{slow}$. We sum the activity of the two oscillators together to produce a multiscale time series of one variable. We adjust the frequency ratio $F=T_\text{slow}/T_\text{fast}$ and compare the performance of the three different HAVOK models discussed above. When adjusting the frequency ratio, we keep the period of the fast oscillator fixed and adjust the period of the slow oscillator. In all models we train using time series data covering five periods of the slow oscillation ($m\Delta t=5T_\text{slow}$). As $F$ increases, we scale our models in three different ways. For the HAVOK model with a large time step, we fix $q$ and increase $\Delta t$ such that the delay duration $D=(q-1)\Delta t=T_\text{slow}$. For the HAVOK model with a small time step, we fix $\Delta t$ and increase the number of delays $q$ to obtain delay duration $D=(q-1)\Delta t=T_\text{slow}$. Finally, for the HAVOK model with spacing, we fix $q,\Delta t$ and increase the row spacing $d$ such that $D=(q-1)d\Delta t=T_\text{slow}$. We also increase column spacing $c$ such that we have the same number of samples for each system. This allows us to maintain a constant size of the snapshot matrix as $F$ increases.

 In the top panel of Figure~\ref{fig:5} we compare both the size of the snapshot matrices and the prediction error for our three different models. By increasing the time step in the standard HAVOK model as the frequency ratio increases, we are able to avoid increasing the size of the snapshot matrix but are penalized with a large increase in error (dashed lines). In contrast, building a standard HAVOK model with a small time step allows us to maintain a small error but causes $F^2$ growth in the size of the snapshot matrices (solid lines). By introducing row and column spacing, we can significantly reduce the error while also avoiding growth in the size of the snapshot matrices (dotted lines). While we choose our row and column spacings $d,c$ to avoid any growth in the size of the snapshot matrices as the frequency ratio grows, these can be adjusted to give more flexibility in the size/error trade off. This allows the method to be very flexible based on computational resources and accuracy requirements.

\subsubsection{Method 2: iterative modeling}
\label{sec:multiscale2_iterative}
Our second strategy is an iterative method for modeling multiscale systems using HAVOK. A key enabling assumption is that with sufficient frequency separation between time scales, the slow dynamics will appear approximately constant in the relevant time scale of the fast dynamics. This allows us to build a model for the fast dynamics using data from a short time period without capturing for the slow dynamics. We then use this model to predict the fast behavior over a longer duration, subtract it out, and model the slow dynamics using downsampled data.

The algorithm proceeds in three steps, illustrated in Figure~\ref{fig:5}. As training data we use the same time series of two Van der Pol oscillators used in Section~\ref{sec:multiscale2_spacing}. In Figure~\ref{fig:5} we show frequency ratio $F=20$. In the first step, we sample for a duration equivalent to five periods of the fast oscillation at a rate of $128$ samples/period. We use this data to build a HAVOK model with rank $r=50$ and $q=128$, using the procedure described in Section~\ref{fig:5a}. The slow dynamics are relatively constant during this recording period, meaning the model only captures the fast dynamics. Setting the real part of all DMD eigenvalues to zero and subtracting any constants (modes with zero eigenvalues), we obtain a model for the fast dynamics.  Such subtraction of {\em background} modes has been used effectively for foreground-background separation in video streams~\cite{jake,Erichson2015arxiv}

The next step is to approximate the slow dynamics, which will be used to build the slow model. We sample the system over a longer duration at a reduced sampling rate; in the example we sample for five periods of the slow oscillation and reduce the sampling rate by a factor of the frequency ratio $F$. We use our fast model to determine the contribution of the fast dynamics at the sampled time points and subtract this out to get an approximation of the slow dynamics. Note that because we have a linear model of the form in (\ref{eq:dmd_model}), we can predict the contribution of the fast oscillator at these sample points without needing to run a simulation on the time scale of the fast dynamics.

Finally we use this subsampled data as training data to construct a HAVOK model of the slow dynamics. In the example in Figure~\ref{fig:5}, we again build a rank 50 model using 128 delay coordinates. This procedure gives us two models, one for the fast dynamics and one for the slow dynamics. These models can be combined to give a model of the full system or used separately to analyze the individual time scales.

This method provides an additional strategy for applying HAVOK to multiscale systems. As it relies on the assumption that the slow dynamics are relatively constant in the relevant time scale of the fast dynamics, it is best suited to problems where there is large separation between time scales. In practice, the application of this method requires a choice of how long to sample to capture the fast dynamics. The method is not sensitive to the exact sampling duration used provided it captures at least a few periods of the fast oscillation and the slow dynamics remain approximately constant throughout this duration. One strategy would be to use the frequency spectrum of the data to inform this decision.

\section{Discussion}
\label{sec:discussion}

The discovery of nonlinear dynamical systems from time series recordings of physical and biological systems has the potential to revolutionize scientific discovery efforts.  Through principled regression methods, the emerging SINDy and HAVOK methods are allowing researchers to discover governing evolution equations by sampling either the full or partial state space of a give system, respectively.   Although these methods are emerging as viable techniques for a broad range of applications, multiscale complex systems remain difficult to characterize for any current mathematical strategy.  Indeed, multiscale modeling is a grand challenge problem for complex systems of the modern area such as climate modeling, brain dynamics, and power grid systems, to name just a few examples of technological importance.

In this manuscript, we have developed a suite of methods capable of extracting interpretable models from multiscale systems.  The methods are effective for building dynamical models when either the full or partial state of the system is measured.  When full state measurements are acquired, then the SINDy sparse regression framework can be modified to discover distinct dynamics at time scales that are well separated in time, i.e. slow and fast dynamics that are coupled.  If only partial measurements of the full state space are recorded, then time-delay embedding and HAVOK can be used to construct a model of the system, despite latent (unknown or unmeasured) state space variables.  As with the multiscale SINDy architecture, the HAVOK architecture can be modified to discover dynamical models at different time scales through the DMD regression framework.

Critical to the success of multiscale discovery are robust temporal sampling strategies of the dynamical system.   We have developed a suite of sampling strategies that allow one to sample either the full or partial state space of the system and discover dynamics in a computationally scalable manner.   Such techniques, in combination with our SINDy and HAVOK multiscale modifications, provide automated strategies for model discovery in complex systems of interest in the modern era.  Moreover, our methods also give guidelines for the minimal amount of data required to produce accurate multiscale models.   For the purposes of reproducibility and to aid in scientific discovery efforts, all code and data has been provided as open source software at github.com/kpchamp/MultiscaleDiscovery.

\section*{Acknowledgments}
JNK acknowledges support from the Air Force Office of Scientific Research (AFOSR) grant FA9550-17-1-0329. SLB acknowledges support from the Army Research Office (W911NF-17-1-0306).  SLB and JNK acknowledge support from the Defense Advanced Research Projects Agency (DARPA contract HR011-16-C-0016). This material is based upon work supported by the National Science Foundation Graduate Research Fellowship under Grant No. DGE-1256082. We are especially thankful to Bing Brunton, Eurika Kaiser and Josh Proctor for conversations related to the SINDy and HAVOK algorithms.


\end{document}